\documentclass{amsart}

\usepackage{graphicx}

\newtheorem{theorem}{Theorem}[section]
\newtheorem{lemma}[theorem]{Lemma}

\theoremstyle{definition}
\newtheorem{definition}[theorem]{Definition}
\newtheorem{example}[theorem]{Example}

\newtheorem{proposition}[theorem]{Proposition}
\newtheorem{corollary}[theorem]{Corollary}

\theoremstyle{remark}
\newtheorem{remark}[theorem]{Remark}

\def\cal{\mathcal}  
\newcommand{\A}{{\cal A}}

\newcommand{\codim}{\mbox{\rm codim}}

\numberwithin{equation}{section}

\def\Der{\mathrm{Der}}

\begin{document}

\title[]{The characteristic polynomial of a multiarrangement}


\author{Takuro Abe}
\address{Department of Mathematics, Hokkaido University, Sapporo, 060-0810, Japan}
\email{abetaku@math.sci.hokudai.ac.jp}

\author{Hiroaki Terao}
\address{Department of Mathematics, Hokkaido University, Sapporo, 060-0810, Japan}
\curraddr{}
\email{terao@math.sci.hokudai.ac.jp}

\author{Max Wakefield}
\address{Department of Mathematics, Hokkaido University, Sapporo, 060-0810, Japan}
\curraddr{}
\email{wakefield@math.sci.hokudai.ac.jp}

\thanks{The first author is been supported by 21st Century COE Program 
``Mathematics of Nonlinear Structures via Singularities'' Hokkaido University. 
The second author
has been
supported in part by Japan Society for the Promotion
of Science. The third author has been supported by NSF grant \# 0600893 and the NSF Japan program. }


\date{\today}

\dedicatory{}

\begin{abstract} Given a multiarrangement of hyperplanes we define a series by sums of the Hilbert series of the derivation modules of the multiarrangement. This series turns out to be a polynomial. Using this polynomial we define the characteristic polynomial of a multiarrangement which generalizes the characteristic polynomial of an arragnement. The characteristic polynomial of an arrangement is a combinatorial invariant, but this generalized characteristic polynomial is not. However, when the multiarrangement is free, we are able to prove the factorization theorem for the characteristic polynomial. The main result is a formula that relates `global' data to `local' data of a multiarrangement given by the coefficients of the respective characteristic polynomials. 
This result gives a new necessary condition 
for a multiarrangement to be free.
Consequently it provides a simple method to show that a given multiarrangement
is not free.
\end{abstract}
\maketitle
\setcounter{section}{-1}
\section{Introduction}
\bigskip
Let $V$ be a vector space of dimension $\ell$ over a field $\mathbb{K}$ and $S=S(V^*)$ be the symmetric algebra. We can choose coordinates for $V^*$ such that $S\cong \mathbb{K}[x_1,\ldots ,x_\ell ]$. Put $\partial_{x_i}:=\partial / \partial x_i$. A hyperplane is a codimension one linear space in $V$. 
A {\it multiarrangement} is a finite collection of hyperplanes 
denoted by $\mathcal{A}$ together with a multiplicity function 
$m:\mathcal{A}\to \mathbb{Z}_{>0}$. 
Let $(\A, m)$ denote a multiarrangement.
When $m(H) = 1$ for all $H\in\A$, we identify $(\A, m)$ with the arrangement
$\A$.    
For $p\geq 1$ the $S$-module $\Der ^p(S)$ is the set of all alternating $p$-linear functions $\theta :S^p\to S$ such that $\theta$ is a $\mathbb{K}$-derivation in each variable. For $p=0$ we put $\Der^0(S)=S$. For each $H\in \mathcal{A}$ we choose a defining form $\alpha_H$. Put $\tilde{Q}=\prod\limits_{H\in \mathcal{A}}\alpha_H^{m(H)}$. Define the derivation modules of $(\mathcal{A},m)$ as 
\begin{multline*}
D^p(\mathcal{A},m)=\{ \theta \in \Der^p(S) | \theta (\alpha_H,f_2,\ldots , f_p)\in \alpha_H^{m(H)}S \\
\text{ for all } H\in \mathcal{A} \text{\ and\ } f_2,\ldots ,f_p\in S
\}.
\end{multline*}
If $D^{1} (\A, m)$ is a free $S$-module we say that
a multiarrangement $(\A, m)$ is {\it free}.  

One of the most fundamental invariants of an arrangement of hyperplanes is its characteristic polynomial. The focus of this paper is to generalize the characteristic polynomial to multiarrangements of hyperplanes and apply this polynomial to the problem of freeness of the module of derivations. In \cite{Z} Ziegler initiated the study of derivations
of multiarrangements.  
Later in \cite{Y2} and \cite{Y3}
Yoshinaga found that 
the derivation modules of multiarrangements
are important for the study of free
arrangements.
It is known that any multiarrangement is free when $\ell = 2$
(see \cite{WY} and \cite{Z}).
Other examples of free multiarrangements include
the restricted multiarrangements of a free arrangement
(see \cite{Z}) and the Coxeter arrangements with a constant multiplicity
(see \cite{T3} and \cite{Y1}).
On the other hand, very few examples of non-free multiarrangements 
have been known.  One purpose of this paper is to introduce a useful criterion
for a multiarrangement to be
non-free. 

In order to define the characteristic polynomial of
a multiarrangement $(\mathcal{A},m)$ we make use of
the $S$-modules $D^p(\mathcal{A},m)$. Since each $D^p(\mathcal{A},m)$ 
is $\mathbb{Z}_{\geq 0}$-graded by polynomial degree,
we may define a function 
$$\psi (\mathcal{A},m;t,q)=\sum\limits_{p=0}^{\ell}H(D^p(\mathcal{A},m),q)(t(q-1)-1)^p$$
in $t$ and $q$, where 
$H(D^p(\mathcal{A},m),q)$
is the Hilbert series of $D^p(\mathcal{A},m)$.
Although
$\psi (\mathcal{A},m;t,q)$
is, a priori,  a rational function
in $q$, it turns out to be a polynomial in $q$
as shown in Theorem \ref{chithm}. 
So we may substitute $q$ equal to
 $1$ and we
define
the {\it characteristic polynomial} 
by
$$\chi ((\mathcal{A},m),t)=(-1)^\ell \psi (\mathcal{A},m;t,1)$$
and the {\it Poincar\'e polynomial} by
$$\pi ((\mathcal{A},m),t)=(-t)^\ell \chi ((\mathcal{A},m),-t^{-1}).$$
These polynomials are generalizations of
the characteristic and Poincar\'e 
polynomials of an arrangement $\A$ because of 
\cite{ST}. 
However, 
unlike the case of arrangements,
these generalized polynomials
are not combinatorial invariants.

Let $L=L(\mathcal{A})$ 
be the intersection lattice of $\mathcal{A}$ with the order 
as reverse inclusion and the rank function defined by 
codimension: $r(X)=\codim_V (X)$. Let $L_k=\{ X\in L
\mid
r(X)=k\}$. 
For any $X\in L$ let 
$\mathcal{A}_X=\{H\in \mathcal{A}
\mid X\subseteq H\}$ 
and $m_X=m|_{\mathcal{A}_X}$. 
Define
$C_p(X)\in \mathbb Z$ 
by
$\pi ((\mathcal{A}_{X} ,m_{X}),t)
=\sum_{p=0}^{\ell}  C_{p}(X) t^{p}$.
The Main Theorem 
\ref{main} in this paper asserts that,
for arbitrary $X\in L$ and $p$ such that 
$0\leq p \leq r(X)$,
$$C_p(X)=\sum\limits_{Y\in L(\A_{X})_{p}}C_p(Y).$$
In particular, when $X$ is equal to the intersection of all
hyperplanes in $\A$,  we have
$$C_p=\sum\limits_{Y\in L_{p}}C_p(Y),$$
where $C_{p}$ is the coefficient of $t^{p}$ 
in the ``global'' Poincar\'e polynomial 
$\pi ((\mathcal{A} ,m),t)$. This formula thus
relates global data to local data 
of derivations of multiarrangements.

The multiset $\exp(\A, m)$ of {\it exponents} are defined by
the polynomial degrees of  a homogeneous 
basis over $S$ as in \cite{Z}
if $(\A, m)$ is a free multiarrangement. 
Next we prove the Factorization Theorem
\ref{factorization}
for free multiarrangements:
$$\pi ((\mathcal{A},m),t)=\prod\limits_{i=1}^\ell (1+d_i t)
$$
where
$\exp (\mathcal{A},m)=(d_1,\ldots ,d_\ell )$.
This is a 
generalization of the factorization theorems
for free arrangements
in
\cite{T1} and \cite{T2}.
When $(\A_{X} , m_{X})$ is free with
$\exp(\A_{X} , m_{X}) =(
d^{X}_{1},
\dots
d^{X}_{r(X)},
0,
\dots,
0)$, 
the Factorization Theorem implies
$C_{r(X)}(X) = 
d^{X}_{1}d_2^X
\cdots
d^{X}_{r(X)}
$.
We define
the $k$-th {\it local mixed product} by
$$LMP(k)=\sum\limits_{X\in L_k}d^{X}_1d^{X}_2\cdots d^{X}_k
$$
when
the multiarrangement $(\mathcal{A}_X,m_X)$ is free
for any $X\in L_k$.
Assuming that $(\A,m)$ is free with
$\exp (\mathcal{A},m)=(d_1,\ldots ,d_\ell )$,
we introduce
the  $k$-th {\it global mixed product} by
$$GMP(k)=\sum d_{i_1}d_{i_2}\cdots d_{i_k}$$ 
where the sum is over all $k$-tuples such that $1\leq i_1<\cdots <i_k\leq \ell$.
Then, thanks to Theorem \ref{main},
we have Corollary \ref{GMP=LMP}:
$$GMP(k)=LMP(k).$$ 
This formula gives a necessary condition 
for a multiarrangement to be free.
Therefore, it provides a simple method to show non-freeness of 
a given multiarrangement as illustrated 
in Example \ref{nonfreeexample}.


\section{Preliminaries}
\bigskip
Let $(\mathcal{A},m)$ be a multiarrangement. 
In this section we collect basic properties of the modules
 $D^p(\mathcal{A},m)$.
We write $(\mathcal{A},m)\subseteq (\mathcal{B},m')$ if $\mathcal{A}\subseteq \mathcal{B}$ and for all $H\in \mathcal{A}\subseteq \mathcal{B}$ we have $0<m(H)\leq m'(H)$. 

\begin{lemma}\label{subset}
If $(\mathcal{A},m)\subseteq (\mathcal{B},m')$ then $D^p(\mathcal{A},m)\supseteq D^p(\mathcal{B},m')$. \end{lemma}

\proof Let $\theta \in D^p(\mathcal{B},m')$ and let $H\in \mathcal{A}$. Then $\theta (\alpha_H,f_2,\ldots ,f_p)\in \alpha_H^{m'(H)}S\subseteq \alpha_H^{m(H)}S$. Thus, $\theta \in D^p(\mathcal{A},m)$. \qed

\
  
We have a product structure on $\Der^p(S)\times \Der^q(S)\to \Der^{p+q}(S)$ because $\bigwedge^p\Der^1(S)\cong \Der^p(S)$. Recall the formula (2.3) from \cite{ST} that if $\theta_1,\ldots ,\theta_p\in \Der^1(S)$ then \begin{equation}\label{det}(\theta_1 \wedge \cdots \wedge \theta_p)(f_1,\ldots,f_p)=\det [\theta_i(f_j)]_{1\leq i,j\leq p}\end{equation} for all $f_1,\ldots ,f_p\in S$. Similarly as described in \cite{ST}, if $\varphi \in D^p(\mathcal{A},m)$ and $\psi \in D^q(\mathcal{A},m)$ then $\varphi \wedge \psi \in D^{p+q}(\mathcal{A},m)$. The next three lemmas are nearly identical to Propositions (2.5), (3.4), and (5.8) in \cite{ST} respectively. However, because they are generalizations and the results are important for this paper we show their proofs.

\begin{lemma}\label{D^l}
$D^\ell (\mathcal{A},m)\cong S\tilde{Q}(\partial_{x_1}\wedge \cdots \wedge \partial_{x_\ell})$.
\end{lemma}

\proof Let $\theta =f(\partial_{x_1}\wedge \cdots \wedge \partial_{x_\ell})\in D^\ell (\mathcal{A},m)$ for some $f\in S$. Let $H\in \mathcal{A}$ be arbitrary. Then for some $i\in \{ 1,\ldots ,\ell \}$ we have that $\partial_{x_i}(\alpha_H)\neq 0$ and without loss of generality we can assume that $i=1$ and $\partial_{x_1}(\alpha_H)=1$ . Then (\ref{det}) implies $\theta (\alpha_H,x_2,\ldots ,x_\ell )=f$. Since $\theta \in D^\ell (\mathcal{A},m)$ we know that $\theta (\alpha_H,x_2,\ldots ,x_\ell )=f\in \alpha_H^{m(H)}S$. So, for all $H\in \mathcal{A}$ we have $f\in \alpha_H^{m(H)}S$. Thus, the polynomial $\tilde{Q}$ divides $f$. \qed

\begin{lemma}\label{wedgee}
If $(\mathcal{A},m)$ is a free multiarrangement then $D^p(\mathcal{A},m)\cong \bigwedge^p D^1(\mathcal{A},m)$.
\end{lemma}

\proof Let $\{ \theta_1,\ldots ,\theta_\ell \}$ be a basis for $D^1(\mathcal{A})$. Let $I=(i_1,\ldots ,i_p)$ where $1\leq i_1<\cdots <i_p\leq \ell$. Let $\partial_I=\partial_{x_{i_1}}\wedge \cdots \wedge \partial_{x_{i_p}}$ and let $\theta_I=\theta_{i_1}\wedge \cdots \wedge \theta_{i_p}$. For all $i\in \{1,\ldots ,\ell\}$ we know $\tilde{Q}\partial_{x_i}\in D^1(\mathcal{A},m)=\sum\limits_{j=1}^\ell S\theta_j$. Thus, $\tilde{Q}^p\Der^p(S)\subseteq \sum\limits_IS\theta_I$. Let $\theta \in D^p(\mathcal{A},m)$. Then there exists $f_I\in S$ such that $\tilde{Q}^p\theta =\sum\limits_If_I\theta_I$. Let $J\subseteq \{ 1,\ldots ,\ell \}$ such that $|J|=\ell -p$. Then by Ziegler's criterion (i.e.,  the multiarrangement version of Saito's criterion, see \cite{Z}) \begin{equation}\label{prodth}\tilde{Q}^p(\theta \wedge \theta_J)=\left( \sum\limits_If_I\theta_I \right) \wedge \theta_J=f_K\tilde{Q}(\partial_{x_1}\wedge \cdots \wedge \partial_{x_\ell})\end{equation} where $K$ is the complement of $J$. Also, $\theta \wedge \theta_J \in D^\ell (\mathcal{A},m)=\tilde{Q}(\partial_{x_1},\ldots ,\partial_{x_\ell})$. Thus, $\tilde{Q}^p$ divides $f_K$ for all $K$. Therefore, $$\theta = \sum\limits_I\frac{f_I}{\tilde{Q}^p}\theta_I.$$ If $\sum\limits_If_I\theta_I=0$ then the second equality of (\ref{prodth}) implies that $f_K=0$ for all $K$. Thus, $\{\theta_I| \ |I|=p\}$ is a basis for $D^p(\mathcal{A},m)$. \qed

\

Let $(\mathcal{A}_1,m_1)$ and $(\mathcal{A}_2,m_2)$ be two multiarrangements in the vector spaces $V_1$ and $V_2$ resprectively.
We define the product of these two multiarrangements by $(\mathcal{A}_1,m_1)\times (\mathcal{A}_2,m_2):=(\mathcal{A}_1\times \mathcal{A}_2,m)$ where the hyperplanes are given by $\mathcal{A}_1\times \mathcal{A}_2=\{ H\oplus V_2| H\in \mathcal{A}_1\}\cup \{ V_1\oplus H'| H'\in \mathcal{A}_2\}\subseteq V_1\oplus V_2$ and the multiplicities are given by $m(H\oplus V_2)=m_1(H)$ and $m(V_2\oplus H')=m_2(H')$. Put $S_i=S(V_i^*)$ for $i=1,2$ and $S=S(V_1^*\oplus V_2^*)$.

\begin{lemma}\label{productss}
$$D^k\left( (\mathcal{A}_1,m_1)\times (\mathcal{A}_2,m_2)\right) \cong \bigoplus\limits_{i+j=k}D^i(\mathcal{A}_1,m_1)\otimes_{\mathbb{K}}D^j(\mathcal{A}_2,m_2).$$
\end{lemma}

\proof In this proof the tensor product is always over $\mathbb{K}$. Identify $S$ with $S_1\otimes S_2$ and $\Der^k(S)$ with $\bigoplus\limits_{i+j=k} \Der^i (S_1)\otimes \Der^j (S_2)$. It is clear that $\bigoplus\limits_{i+j=k}D^i(\mathcal{A}_1,m_1)\otimes D^j(\mathcal{A}_2,m_2) \subseteq D^k\left( (\mathcal{A}_1,m_1)\times (\mathcal{A}_2,m_2)\right) $. We show the reverse inclusion. Let $\theta \in D^k((\mathcal{A}_1,m_1)\times (\mathcal{A}_2,m_2))$. Without loss of generality we can assume that $\theta \in \Der^i (S_1)\otimes \Der^j (S_2)$. Suppose that $\theta =\sum\limits_{s=1}^r\varphi_s\otimes \psi_s$ where $\varphi_1,\ldots ,\varphi_r\in \Der^i(S_1)$ and $\psi_1,\ldots , \psi_r\in \Der^j(S_2)$ are linearly independent over $\mathbb{K}$. Fix $f_2,\ldots , f_i\in S_1$ and for any $H\in \mathcal{A}_1$ let $\psi = \sum\limits_{s=1}^r\varphi_s(\alpha_H,f_2,\ldots ,f_i)\otimes \psi_s \in S_1\otimes \Der^j(S_2)\subseteq \Der^j(S)$. Let $(\Phi ,m_{\emptyset})$ be the empty multiarrangenment in $V_2$. Since $\theta \in D^k((\mathcal{A}_1,m_1)\times (\mathcal{A}_2,m_2))\subseteq D^k((\mathcal{A}_1,m_1)\times (\Phi ,m_{\emptyset}))$ we know that for all $g_1,\ldots , g_j\in S_2$ we have $$\psi (g_1,\ldots ,g_j )=\left( \sum\limits_{s=1}^r\varphi_s\otimes \psi_s\right)(\alpha_H,f_2,\ldots ,f_i,g_1,\ldots ,g_j)\in \alpha_H^{m_1(H)}S.$$ Thus, $\psi  \in \alpha_H^{m_1(H)}S_1\otimes \Der^j(S_2)$. Since $\psi_1,\ldots ,\psi_r$ are linearly independent over $\mathbb{K}$ we know that $1\otimes \psi_1,\ldots ,1\otimes \psi_r$ are linearly independent over $S_1$. Therefore, $\varphi_s(\alpha_H,f_2,\ldots,f_i)\in \alpha_H^{m_1(H)}S_1$ and $\theta \in D^i(\mathcal{A}_1,m_1)\otimes \Der^j(S_2)$. Now, we can choose $\xi_1,\ldots ,\xi_t\in D^i(\mathcal{A}_1,m_1)$ that are linearly independent over $\mathbb{K}$ such that $\theta =\sum\limits_{s=1}^t\xi_s\otimes \zeta_s$ for some $\zeta_1,\ldots ,\zeta_t \in \Der^j(S_2)$. To finish the proof we just perform the same argument to the $\zeta_1,\ldots ,\zeta_t$ as we did above with the $\varphi_s$ and we have that $\theta \in D^i(\mathcal{A}_1,m_1)\otimes D^j(\mathcal{A}_2,m_2)$. \qed

\

Let $(S\mathrm{-Mod})$ denote the category of $S$-modules. Regard $L$ as a category with morphisms $\leq$. Next we follow \cite{RT} and using the modules $D^p(\mathcal{A},m)$ we define a contravariant functor $$D^p:L\to (S\mathrm{-Mod})$$ by $D^p(X):=D^p(\mathcal{A}_X,m_X)$ and $D^p(\leq )$ is the inclusion from Lemma \ref{subset}. We review the definition of a local functor from \cite{ST}.

\begin{definition}
For any prime ideal $P\subseteq S$ let $X(P)=\bigcap H$ where the intersection is over all $H\in \mathcal{A}$ such that $X\subseteq H$ and $\alpha_H\in P$. We say that a contravariant functor $F:L\to (S\mathrm{-Mod})$ is {\it local} if the localization of $F(X)\to F(X(P))$ at $P$ is an isomorphism for every $X\in L$ and every prime ideal $P$.
\end{definition}

Now the proof that $D^p$ is a local functor is slightly different from the proof in \cite{RT}.

\begin{proposition}\label{Dlocal}
For every $0\leq p\leq \ell$ the functors $D^p$ are local.
\end{proposition}

\proof Let $P$ be a prime ideal of $S$. For every $0\leq p\leq \ell$ we have the inclusion $D^p(\mathcal{A}_{X(P)},m_{X(P)})\supseteq D^p(\mathcal{A}_X,m_X)$ by Lemma \ref{subset}. Let $$\frac{\theta}{f}\in D^p(\mathcal{A}_{X(P)},m_{X(P)})_P$$ where $\theta \in D^p(\mathcal{A}_{X(P)},m_{X(P)})$ and $f\in S\backslash P$. Define the polynomial $$g=\prod\limits_{H\in \mathcal{A}_X\backslash \mathcal{A}_{X(P)}}\alpha_H^{m(H)}.$$ Then $$\frac{g\theta}{gf}\in D^p(\mathcal{A}_X,m_X)_P.$$ Thus, $D^p(\mathcal{A}_X,m_X)_P\cong D^p(\mathcal{A}_{X(P)},m_{X(P)})_P.$ \qed

\

The following theorem from \cite{ST} is crucial in the proof of our main result. To state it we need to have some notation. Each $D^p(\mathcal{A},m)$ is $\mathbb{Z}_{\geq 0}$-graded by the polynomial grading. The Hilbert series (also called the Poincar\'e series) of an $\mathbb{Z}_{\geq 0}$-graded, finitely generated module $M=\bigoplus M_q$ is $$H(M,q)=\sum\limits_{p=0}^{\infty}\dim_\mathbb{K}(M_p)q^p.$$
Let $\mu : L\times L \longrightarrow \mathbb Z$ be the M\"obius function
as in \cite{Rota} and
\cite{Stanley}.

\begin{theorem}[\cite{ST}, (6.10)]\label{polethm}
Let $F$ be a contravariant, $\mathbb{Z}_{\geq 0}$-graded, finitely generated, local functor $F:L \to (S\mathrm{-Mod})$. Then for any $X\in L$ $$\sum\limits_{Y\leq X}\mu (Y,X)H(F(Y),q)$$ has a pole of order at most $\dim X$ at $q=1$. 
\end{theorem}


\section{Definition of $\chi ((\mathcal{A},m ),t)$}\label{chii}

Let $(\mathcal{A},m)$ be any multiarrangement. In this section we define a series $\psi (\mathcal{A},m;t,q)$ associated to the multiarrangement $(\mathcal{A},m)$, prove that it is a polynomial, and then with this polynomial define the characteristic polynomial $\chi ((\mathcal{A},m),t)$ and the Poincar\'e polynomial $\pi ((\mathcal{A},m),t)$.

\begin{definition}\label{generalpsi}
$$\psi (\mathcal{A},m;t,q)=\sum\limits_{p=0}^{\ell}H(D^p(\mathcal{A},m),q)(t(q-1)-1)^p.$$
\end{definition}

Next we summarize the arguments in \cite{ST} for the case of multiarrangements to prove that $\psi (\mathcal{A},m;t,q)$ is a polynomial in $q$ and $t$. The symmetric algebra $S(V^*)$ is $\mathbb{Z}_{\geq 0}$-graded by homogeneous polynomial degree and we denote the $d$-th graded component by $S(V^*)_d$.

\begin{definition}\label{nondegenerate}
We say $h\in S(X^*)_d$ is {\it non-degenerate} on a subspace $X\subseteq V$ if $\sqrt{(\partial_{x_1}(h),\ldots ,\partial_{x_k}(h))}\supseteq (x_1\ldots ,x_k)$ where $\{ x_1,\ldots ,x_k\}$ is a basis for $X^*$. Let $N^X_d$ be the set of all $h\in S(V^*)_d$ such that  $h|_X\in S(X^*)_d$ is non-degenerate on $X\subseteq V$.
\end{definition}

\begin{remark}
Assume $\mathbb{K}$ is algebraically closed, then there are infinitely many $d$ such that $N^X_d$ is non-empty and actually a Zariski open set in $S(X^*)_d$ for all $X\in L$ (see \cite{ST}).

\end{remark}

\begin{lemma}\label{radical}
Assuming $\mathbb{K}$ is algebraically closed, if  $h\in \bigcap\limits_{X\in L}N^X_d$ then the ideal $\sqrt{D(\mathcal{A},m)h}\subseteq S(V^*)$ contains the unique homogeneous maximal ideal.
\end{lemma}

\proof By Hilbert's Nullstellensatz it is enough to show that the zero locus $V(D(\mathcal{A},m)h)$ is contained in $\{ 0\}$. For $v\in V\backslash \{ 0\}$ let $X=\bigcap\limits_{v\in H\in \mathcal{A}}H$. So, $v\in X$, but $v\notin Y$ for all $Y\in L$ such that $Y\subset X$. Choose a basis $\{ x_1,\ldots ,x_\ell\}$ for $V^*$ such that $X=V(x_{k+1},\ldots ,x_\ell )$. Let $\tilde{Q}'=\prod\limits_{H\not\supseteq X}\alpha_H^{m(H)}$. It is clear that $\tilde{Q}'\partial_{x_i}\in D^1(\mathcal{A},m)$ for all $i\in \{ 1,\ldots ,k\}$. Since $h\in N^X_d$ we know that $v\notin V(\partial_{x_1} (h),\ldots ,\partial_{x_k} (h))\cap X$ so that $v\notin V(\partial_{x_1} (h),\ldots ,\partial_{x_k} (h))$. Since $X\in L$ is minimal such that $v\in X$ we have that $\tilde{Q}'(v)\neq 0$. Thus, $v\notin V(\tilde{Q}'\partial_{x_1}(h),\ldots ,\tilde{Q}'\partial_{x_k}(h))$. But $V(D(\mathcal{A},m)h)\subseteq V(\tilde{Q}'\partial_{x_1}(h),\ldots ,\tilde{Q}'\partial_{x_k}(h))$ so $v\notin V(D(\mathcal{A},m)h)$. \qed

\begin{theorem}\label{chithm}
The series $\psi (\mathcal{A},m;t,q)$ is a polynomial in $q$ and $t$.
\end{theorem}

\proof First, we note that since $\psi (\mathcal{A},m;t,q)$ is stable under field extension, we may assume that $\mathbb{K}$ is algebraically closed. Consider the following chain complex \begin{equation}0\to D^\ell (\mathcal{A},m)\to D^{\ell -1}(\mathcal{A},m)\to \cdots \to D^1(\mathcal{A},m)\to D^0(\mathcal{A},m)\to 0 \end{equation} where the differential $\partial_h$ is defined by $$(\partial_h \theta )(f_1,\ldots ,f_{p-1}):=\theta (h,f_1,\ldots ,f_{p-1})$$ for any $\theta \in D^p(\mathcal{A},m)$. Once we replace the above complex for the corresponding complex defined in equation (4.7) in \cite{ST} the proof follows from the above Lemma \ref{radical} and Propositions (4.10), (5.2), and (5.3) in \cite{ST}. \qed

\

Because $\psi (\mathcal{A},m;t,q)$ is a polynomial for any multiarrangement we can make the following definition.

\begin{definition}\label{chidef}
The {\it characteristic polynomial} of any multiarrangement $(\mathcal{A},m)$ is the polynomial $$\chi ((\mathcal{A},m),t)=(-1)^\ell \psi (\mathcal{A},m;t,1)$$  and the {\it Poincar\'e polynomial} is $$\pi ((\mathcal{A},m),t)=(-t)^\ell \chi ((\mathcal{A},m),-t^{-1}).$$
\end{definition}

\begin{remark} This generalizes the characteristic and Poincar\'e polynomials of an arrangement because of Theorem (1.2) of \cite{ST}. However, this polynomial $\chi ((\mathcal{A},m),t)$ is in no way ``characteristic'' since it is not an invariant of the intersection lattice (see the next example). 
\end{remark}

\begin{example}[Ziegler \cite{Z}]

Let $(\mathcal{A}_1,m_1)$ and $(\mathcal{A}_2,m_2)$ be defined by the polynomials $\tilde{Q}_1=x^3y^3(x-y)(x+y)$  and $\tilde{Q}_2=x^3y^3(x-y)(x-cy)$ where $c\neq \infty , 0, 1$. Then the characteristic polynomials are: $$\chi ((\mathcal{A}_1,m_1),t)=(t-3)(t-5)$$ and $$\chi ((\mathcal{A}_2,m_2),t)=(t-4)^2.$$

\end{example}


\section{Local to Global formula for $\chi ((\mathcal{A},m),t)$}

By Theorem \ref{chithm} we know that $\psi (\mathcal{A},m;t,q)$ is a polynomial for any multiarrangement. Thus, the following functions are well-defined.

\begin{definition}\label{localfuncs}
For $p\in \{ 0,\ldots ,\ell \}$ define the functions $C_p:L (\mathcal{A})\to \mathbb{Z}$ by setting $C_p(X)$ equal to the coefficient of $t^p$ in the polynomial $\pi ((\mathcal{A}_{X} , m_{X} ),t)$ or equivalently the coefficient of $t^{\ell -p}$ in the polynomial $(-1)^\ell \chi ((\mathcal{A}_X,m_X),-t)$.
\end{definition}

With this notation $\psi (\mathcal{A}_X,m_X;-t,1)=\sum\limits_{p=0}^\ell C_p(X)t^{\ell -p}$ for all $X\in L$. Let $(\Phi_n,m_\emptyset)$ be the empty multiarrangement in dimension $n$.

\begin{remark}
By Lemma \ref{productss} and the fact that $\chi ((\Phi_n,m_\emptyset),t)=t^n$ (see \cite{OT}) we know that $\chi ((\mathcal{A}_X,m_X),t)$ is divisible by $t^{\dim X}$. Therefore, $C_p(X)=0$ for all $p$ such that $\ell-p<\dim X$. 
\end{remark}

Now, we can state the main theorem which, simply put, states that there is a direct relationship between the local data and the global data of derivations on multiarrangements.

\begin{theorem}\label{main}

For arbitrary $X\in L$ and $p$ such that $0\leq p \leq r(X)$  $$C_p(X)=\sum\limits_{Y\in L(\mathcal{A}_X)_{p}}C_p(Y).$$
\end{theorem}

\proof Let $$\psi_X (\mathcal{A},m;t,q):= \sum\limits_{p=0}^\ell \sum\limits_{Y\leq X} \mu (Y,X) H(D^p(\mathcal{A}_Y,m_Y),q)(t(1-q)-1)^p.$$ By interchanging sums and using Definition \ref{generalpsi} \begin{equation}\psi_X(\mathcal{A},m;t,q)=\sum\limits_{Y\leq X}\mu (Y,X)\psi (\mathcal{A}_X,m_X;-t,q).\end{equation} Thus, by setting $q=1$ and using Definition \ref{localfuncs} \begin{equation}\label{cp's}\psi_X(\mathcal{A},m;t,1)=\sum\limits_{Y\leq X}\mu (Y,X)\sum\limits_{p=0}^\ell C_p(Y)t^{\ell -p}.\end{equation}

Examine the following series $$M_p(q):=\sum\limits_{Y\leq X}\mu (Y,X)H(D^p(\mathcal{A}_Y,m_Y),q).$$ By Proposition \ref{Dlocal} and Theorem \ref{polethm} we have that $(1-q)^{\dim X}M_p(q)$ does not have a pole at $q=1$. So, the coefficient of $t^n$ in $M_p(q)(t(1-q)-1)^p$ is divisible by $(1-q)$ for $n>\dim X$. Hence, the coefficient of $t^n$ in $\psi_X (\mathcal{A},m;t,1)$ is zero for $n>\dim X$. 

On the other hand, $\psi_{X} (\A, m; t, 1)$
is divisible by $t^{\dim X} $ because
$C_{p}(Y) = 0 $ in (\ref{cp's}) 
when 
$\ell - p < \dim Y$.
Thus
$
\psi_{X}(\A, m; t, 1) 
$ 
is a monomial of degree $\dim X$.
Comparing the coefficients of $t^{\dim X} $ in
both sides of
 (\ref{cp's}),
we obtain 
\[
\psi_{X}(\A, m; t, 1) 
=
C_{r(X)} (X) t^{\dim X}. 
\]
 The M\"obius inversion formula 
converts 
\[
\sum_{Y\leq X} 
\mu(Y, X)
\psi(\A_{Y}, m_{Y}; -t, 1)
=
C_{r(X)} (X) t^{\dim X} 
\]
into
\[
\psi(\A_{X}, m_{X}; -t, 1)
=
\sum_{Y\leq X} C_{r(Y)} (Y) t^{\dim Y}. 
\]
This completes the proof since
$
\psi(\A_{X}, m_{X}; -t, 1)
=
\sum_{p=0}^{\ell} C_{p} (X) t^{\ell-p}
$. \qed


\section{$\chi ((\mathcal{A},m ),t)$ for free multiarrangements}\label{freee}

In this section we study methods of applying $\chi ((\mathcal{A},m ),t)$ to the problem of determining the freeness of multiarrangements. First, we prove the ``Factorization Theorem'' for multiarrangements (this generalizes the main Theorems of \cite{T1} and \cite{T2}).

\begin{theorem}\label{factorization}
If $D^1(\mathcal{A},m)$ is free with exponents $\exp (\mathcal{A},m)=(d_1,\ldots ,d_\ell )$ then $$\chi ((\mathcal{A},m),t)=\prod\limits_{i=1}^\ell (t-d_i)$$ and $$\pi ((\mathcal{A},m),t)=\prod_{i=1}^\ell (1+d_it).$$
\end{theorem}

\proof Since the module $D^1(\mathcal{A},m)$ is free, we apply Lemma \ref{wedgee} and get that $$H(D^p(\mathcal{A},m),q)=\sum\frac{q^{d_{i_1}+ d_{i_2} +\cdots +d_{i_p}}}{(1-q)^\ell }$$ where the sum is over all $p$-tuples such that $1\leq i_1 < i_2<\cdots <i_p\leq \ell$. Then multiplying by $t(q-1)-1$ in $\psi (\mathcal{A},m;t,q)$ we factor to get that $$\psi (\mathcal{A},m;t,q)=\frac{\prod\limits_{i=1}^\ell (1+q^{d_i}(t(q-1)-1))}{(1-q)^\ell}.$$ Then expand this and divide by $(1-q)$ in each factor to get $$\psi (\mathcal{A},m;t,q)=\prod\limits_{i=1}^\ell (1+q+q^2+\cdots +q^{d_i-1}-q^{d_i}t).$$ Now, we substitute $q=1$ to get \begin{equation*}\label{q=1} \psi (\mathcal{A},m;t,1)= \prod\limits_{i=1}^\ell (d_i-t) . \qedhere \end{equation*}

Now, we construct a different formula for $\chi ((\mathcal{A},m),t)$. Let $X\in L$. Suppose that $(\mathcal{A}_X,m_X)$ is free with exponents  $$\exp (\mathcal{A}_X,m_X)=(d^X_1,\ldots ,d^X_\ell ).$$ Some of the $d^X_i$ may be zero because of Lemma \ref{productss}. Without loss of generality we may assume that $d^X_k=0$ for all $k>r(X)$. Then by Theorem \ref{factorization} $\pi ((\mathcal{A}_X,m_X),t)=\prod\limits_{i=1}^{r(X)}(1+d_i^Xt)$. Applying this to Definition \ref{localfuncs} we have \begin{equation}\label{Cform}C_{r(X)}(X)=d_1^X\cdots d_{r(X)}^X.\end{equation}

The next definition is a generalization of the ideas of locally free arrangements in \cite{MS} and \cite{Y}.

\begin{definition}\label{locallyfree}
Let $0\leq p\leq \ell$. We say $(\mathcal{A},m)$ is $p$-{\it locally free} if for all $0\leq k \leq p$ and for any $X\in L_k$ the multiarrangement $(\mathcal{A}_X,m_X)$ is free.  
\end{definition}

\begin{definition}\label{LMP}
Suppose that $(\mathcal{A},m)$ is $p$-locally free and that $0\leq k\leq p$. The $k$-th {\it local mixed product} is $$LMP(k)=\sum\limits_{X\in L_k}d^{X}_1d^{X}_2\cdots d^{X}_k.$$
\end{definition}

Since every multiarrangement is $2$-locally free $LMP(2)$ is always well-defined. The next corollary directly follows from (\ref{Cform}), Definitions \ref{locallyfree} and \ref{LMP}, and Theorem \ref{main}.

\begin{corollary}\label{localform}
If $(\mathcal{A},m)$ is $p$-locally free then for all $0\leq k\leq p$ the coefficient of $t^k$ in $\pi ((\mathcal{A},m),t)$ is $LMP(k)$.
\end{corollary}

\begin{definition}\label{GMP}
Let $0\leq k\leq\ell$ and let $(\mathcal{A},m)$ be a free multiarrangement with $\exp (\mathcal{A},m)=(d_1,\ldots ,d_\ell )$. The  $k$-th {\it global mixed product} is $$GMP(k)=\sum d_{i_1}d_{i_2}\cdots d_{i_k}$$ where the sum is over all $k$-tuples such that $1\leq i_1<\cdots <i_k\leq \ell$.
\end{definition} 

Then applying Theorem \ref{factorization} and Corollary \ref{localform} to a free multiarrangement in the setting of Definitions \ref{LMP} and \ref{GMP} we get the following corollary.

\begin{corollary}\label{GMP=LMP}
If $(\mathcal{A},m)$ is a free multiarrangement with $\exp (\mathcal{A},m)=(d_1,\ldots ,d_\ell )$ then for all $0\leq k\leq \ell$ $$GMP(k)=LMP(k).$$ 
\end{corollary}

Now, we describe a simple method to show non-freeness of some multiarrangements. Let $(\mathcal{A},m)$ be a free multiarrangement with exponents $(b_1,\ldots ,b_\ell )$ where $b_1\leq \cdots \leq b_\ell$. Suppose $(d_1,\ldots ,d_\ell )$ is a set of integers such that $d_1\leq \cdots \leq d_\ell$, $\sum\limits_{i=1}^\ell d_i=\sum\limits_{i=1}^\ell b_i$, and $\sum\limits_{i=1}^{\ell -1} (d_{i+1}-d_i) \leq \sum\limits_{i=1}^{\ell -1} (b_{i+1}-b_i)$. We say $(d_1,\ldots ,d_\ell )$ is  ``more balanced'' than $(b_1,\ldots ,b_\ell )$. It is easy to see that $\sum\limits_{i=1}^k d_{i_1}\cdots d_{i_k}\geq \sum\limits_{i=1}^k b_{i_1}\cdots b_{i_k}=GMP(k)$.  By Corollary \ref{GMP=LMP} if $LMP(k)>\sum\limits_{i=1}^\ell d_{i_1}\cdots d_{i_p}\geq GMP(k)$ then we have a contradiction and $(\mathcal{A},m)$ cannot be free.


\section{Applications and Examples}

In \cite{Z} Ziegler shows that there exists a free arrangement $\mathcal{A}$ such that the multiarrangement $(\mathcal{A},m)$ is not free for some multiplicity function $m$ (Example 14). Ziegler proved this by calculating generators for $D^1(\mathcal{A},m)$. Since this is a complicated calculation he asked in \cite{Z} if there is a more systematic method to prove the multiarrangement is not free. Section \ref{freee} provides such a method and we exhibit 
the method on Ziegler's Example 14 in \cite{Z} below.

\begin{example}[Ziegler \cite{Z}]
\label{nonfreeexample}

Let $\mathcal{A}$ be an arrangement defined by the polynomial $Q=xy(x-y)(x-z)(y-z)$. Then $\mathcal{A}$ is free. Let $(\mathcal{A},m)$ be multiarrangement defined by the polynomial $\tilde{Q}=Q^2$. Notice that $(3,3,4)$ is ``more balanced'' than the exponents of $(\mathcal{A},m)$. Thus, in this case $GMP(2)\leq  
3\times 3 
+
3\times 4 
+
3\times 4 
=
33$. 

Figure \ref{maybenotuseit} is the projectivized picture of $(\mathcal{A},m)$ where the circled numbers are the product of the exponents at the corresponding rank two lattice element and the non-circled numbers are the multiplicity of the corresponding line. Summing the products of the exponents at the rank two lattice elements we get that $LMP(2)=34>33\geq GMP(2)$. Therefore, $(\mathcal{A},m)$ is not free.

\renewcommand{\figurename}{Figure  5.\!\!}
\begin{figure}[htbp]\label{maybenotuseit}
\begin{center}
\resizebox{2.7in}{!}{\includegraphics{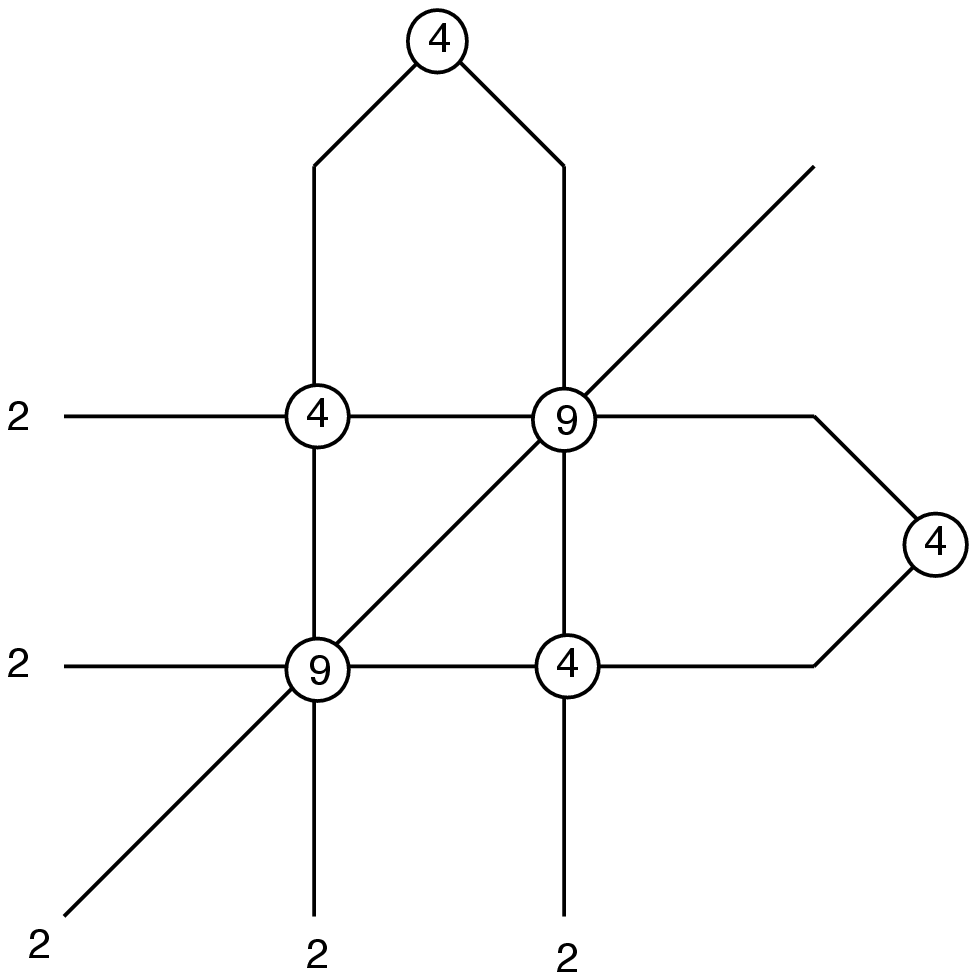}}
\caption{}
\end{center}
\end{figure}

Also, we show a submultiarrangement of the above $(\mathcal{A},m)$ is not free. Let $(\mathcal{A},m_1)$ be a multiarrangement defined by the polynomial $\tilde{Q}_1=x^2y(x-y)(x-z)(y-z)^2$. Suppose $(\mathcal{A},m_1)$ is a free multiarrangement. Then $(2,2,3)$ is ``more balanced'' than the exponents of $(\mathcal{A},m_1)$. Thus, in this case, $GMP(2)
\leq
2\times 2 
+
2\times 3 
+
2\times 3 
=
16$. However, $LMP(2)=17$. Thus, $(\mathcal{A},m_1)$ is not free.
\end{example}

\

One of the most useful theorems concerning the characteristic polynomial of an arbitrary arrangement is the ``Deletion-Restriction'' Theorem (Theorem 2.56 and Corollary 2.57 of \cite{OT}). The theorem states that for any triple $(\mathcal{A},\mathcal{A}',\mathcal{A}'')$ the characteristic polynomials satisfy $\chi (\mathcal{A},t)=\chi (\mathcal{A}',t)-\chi (\mathcal{A}'',t)$. The only suitable generalization of this ``Deletion-Restriction'' theorem to multiarrangements, because of the multiarrangement version of the  `Addition-Deletion'  Theorem in \cite{ATW}, would be that $\chi ((\mathcal{A},m),t)=\chi ((\mathcal{A}',m'),t)-\chi ((\mathcal{A}'',m^*),t)$ where $(\mathcal{A}'',m^*)$ is defined by the `$e$-multiplicities' given in \cite{ATW}. In the next example we show that this generalized `Deletion-Restriction' theorem does not hold for all multiarrangements. Also, this example shows that the characteristic polynomial of a multiarrangement does not necessarily have a linear factor with integer coefficients as is the case for arrangements.

\begin{example}

Let $(\mathcal{A},m)$ be defined by the polynomial $\tilde{Q}=x^2y^2z(x+y+z)(x-y+z)$. Also, let $H_0=\{ y=0\}$ so that $(\mathcal{A}',m')$ is defined by the polynomial $\tilde{Q}'=x^2yz(x+y+z)(x-y+z)$. Figure \ref{maybenotuseit2} is a projectivized picture of $(\mathcal{A},m)$ where the outer circle is the hyperplane at infinity, the un-boxed numbers are the multiplicities of the corresponding projective line and the boxed numbers on $H_0$ are the `$e$-multiplicities' of the corresponding point in the restricted multiarrangement $(\mathcal{A}'',m^*)$. 

\renewcommand{\figurename}{Figure  5.\!\!}
\begin{figure}[htbp]\label{maybenotuseit2}
\begin{center}
\resizebox{3in}{!}{\includegraphics{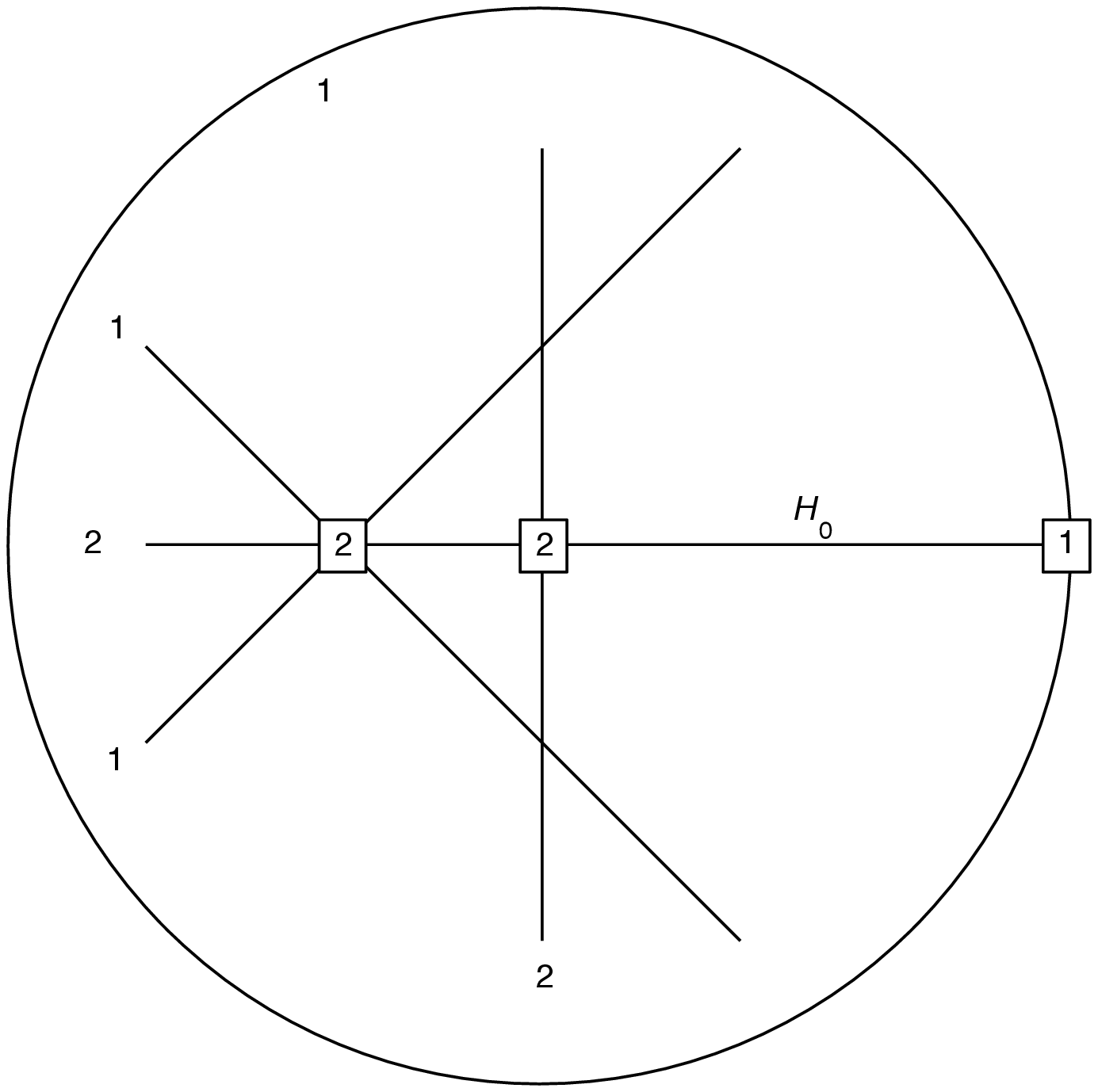}}
\caption{}
\end{center}
\end{figure}

In this case the restricted multiarrangement has $\chi ((\mathcal{A}'',m^*),t)=(t-2)(t-3)$. Then using the computer algebra system {\it Macaulay 2} $$ \chi ((\mathcal{A}',m'),t)=(t-2)(t^2-4t+5)$$ and $$\chi ((\mathcal{A},m),t)=t^3-7t^2+18t-17.$$ Thus, the generalized `Deletion-Restriction' formula is not true for this example and $ \chi ((\mathcal{A},m),t)$ has no integer factor.

\end{example}




\bigskip

\bibliographystyle{amsplain}

\end{document}